\documentclass[12pt,a4paper]{article}
\usepackage[latin2]{inputenc}
\usepackage{amsmath}
\usepackage{amsfonts}
\usepackage{amssymb}
\usepackage{graphicx}
\usepackage[left=2.00cm, right=2.00cm, top=2.50cm, bottom=2.50cm]{geometry}
\usepackage[T1]{fontenc}

\def\Q{{\mathbb Q}}
\def\Z{{\mathbb Z}}

\newtheorem{lemma}{Lemma}
\newtheorem{theorem}[lemma]{Theorem}

\title{
Investigating monogenity in a family of cyclic sextic fields
}
\author{
Istv\'an Ga\'al\\
{\small University of Debrecen, Mathematical Institute} \\
{\small H--4002 Debrecen Pf.400., Hungary,} \\
{\small e--mail: gaal.istvan@unideb.hu},
}

\begin{document}
\baselineskip=17pt

\maketitle
\thispagestyle{empty}

\renewcommand{\thefootnote}{\arabic{footnote}}
\setcounter{footnote}{0}

\vspace{0.5cm}


\begin{abstract}
Jones \cite{jones}  characterized among others monogenity of a family 
of cyclic sextic polynomials. Our purpose is to study monogenity of the family
of corresponding sextic number fields. This also provides the first 
non-trivial application of the method described in \cite{totrelext},
emphasizing its efficiency. 
\end{abstract}

\section{Introduction}

Monogenity and power integral bases is a classical topic in algebraic number theory,
which is intensively studied even today. For the classical results 
we refer to \cite{nark}, for several recent results to \cite{book}. Even during the
last couple of years a huge number of results contributed this area \cite{axioms}.
Here we shorly summarize the most important definitions.

A number field $K$ of degree $n$ with ring of integers $\Z_K$ 
is called {\it monogenic} (cf. \cite{book}) if there 
exists $\xi\in \Z_K$ such that $(1,\xi,\ldots,\xi^{n-1})$ is an integral basis, 
called {\it power integral basis}. 
We call $\xi$ the {\it generator} of this power integral basis.
$\alpha,\beta\in\Z_K$ are called {\it equivalent},
if $\alpha+\beta\in\Z$ or $\alpha-\beta\in\Z_K$. 
Obviously, $\alpha$ generates a power integral basis in $K$
if and only if any $\beta$, equivalent to $\alpha$, does.
As it is known, any algebraic number field admits, up to equivalence, only
finitely many generators of power integral bases.

An irreducible polynomial $f(x)\in\Z[x]$ is called {\it monogenic}, 
if a root $\xi$ of $f(x)$ generates a power integral basis in $K=\Q(\xi)$.
If $f(x)$ is monogenic, then $K$ is also monogenic, but the converse is
not true.

For $\alpha\in\Z_K$ (generating $K$ over $\Q$)  the module index 
\[
I(\alpha)=(\Z_K:\Z[\alpha])
\]
is called the {\it index} of $\alpha$. The element 
$\alpha$ generates a power integral basis in $K$
if and only if $I(\alpha)=1$. If $\alpha^{(i)}\; (1\le i\le n)$ 
are the conjugates of 
$\alpha$ in $K$  of degree $n$, then
\[
I(\alpha)=\frac{1}{\sqrt{|D_K|}}\prod_{1\le i<j\le n}|\alpha^{(i)}-\alpha^{(j)}|,
\]
where $D_K$ is the discriminant of $K$.
Searching for elements of $\Z_K$, generating power integral bases, leads to
a Diophantine equation. Let $(1,\omega_2,\ldots,\omega_n)$ be an integral basis
of $K$, and let 
$L^{(i)}(x_2,\ldots,x_n)=x_2\omega_2^{(i)}+\ldots x_n\omega_n^{(i)}$ 
be the conjugates of the
linear form $L(x_2,\ldots,x_n)=x_2\omega_2+\ldots x_n\omega_n$ ($1\le i\le n$).
The polynomial
\[
I(x_2,\ldots,x_n)=\frac{1}{\sqrt{|D_K|}}
\prod_{1\le i<j\le n}(L^{(i)}(x_2,\ldots,x_n)-L^{(j)}(x_2,\ldots,x_n))
\]
has rational integer coefficients and degree $n(n-1)/2$. It is
called the {\it index form} 
corresponding to the integral basis $(1,\omega_2,\ldots,\omega_n)$.
Obviously, $\gamma=x_1+\omega_2x_2+\ldots+\omega_nx_n\in\Z_K$
generates a power integral basis in $K$ if and only if $(x_2,\ldots,x_n)\in\Z^{n-1}$
is a solution of the {\it index form equation}
\[
I(x_2,\ldots,x_n)=\pm 1
\]
(independently of $x_1$).

Recently several authors apply the method of Newton polygons and Dedeking criterion
which can be used to prove monogenity or non-monogenity of polynomials and number fields.
It is also an important problem to determine explicitly all inequivalent generators 
of power integral bases of a number field. 
The resolution of the index form equations requires numerical algorithms.
Up to now, we only have general efficient algorithms for determining
all inequivalent generators of power integral bases 
only in cubic, quartic and special types of 
sextic and octic number fields. The reason is, that in these cases the relevant
index form equations lead to Thue equations or relative Thue equations, 
which can be easily solved (cf. \cite{book}), at least when their 
coefficients are moderate.

These algorithms for the "complete resolution" of index form equations,
may require too long CPU time. There are also some very fast methods for 
determining generators of power integral bases with "small" coefficients,
say, being $<10^{100}$ in absolute value, with respect to an integral basis.
These solutions cover all solutions with high probability, 
certainly all generators that
can be used in practice for further calculations. 
It is usual to apply such algorithms also if we need to solve a large number of 
equation.

The general method for the complete resolutions of index form equations
in sextic fields
\cite{book} requires a huge amount of time. However, for some special types
of sextic fields we have developed very efficient methods. One of the most 
interesting case is represented by composites of real cubic and imaginary quadratic fields
\cite{totrelext}. In this paper we give the first non-trivial 
application of that method
to a parametric family of cyclic sextic fields, studied by L. Jones \cite{jones}.

For further results on the monogenity of cyclic sextic fields we refer to
\cite{2000}, \cite{2018}, \cite{2019}.

\section{The family of cyclic sextic fields}

Let $n\in\Z$ and consider
\begin{equation}
f(x)=x^6+(n^2+5)x^4+(n^2+2n+6)x^2+1.
\label{f}
\end{equation}

Jones \cite{jones} proved:
\begin{lemma}
$f(x)$ is irreducible for all $n\in\Z$, and is monogenic exactly for $n=-2,-1,0,1$.
\end{lemma}

Let $\xi$ be a root of $f(x)$ and set $K=\Q(\xi)$ with ring of integers $\Z_K$ and 
discriminant $D_K$.
We prove that in case $f(x)$ is monogenic, often 
there are further inequivalent generators of
power integral bases of $K$ in addition to $\xi$. Moreover,
there are parameters $n$, for which $f(x)$ is not monogenic, but $K$ is monogenic.

Let 
\begin{equation}
g(x)=x^3+(n^2+5)x^2+(n^2+2n+6)x+1.
\label{g}
\end{equation}
Denote by $\alpha$ a root of $g(x)$, let $L=\Q(\alpha)$ with ring of integers $\Z_L$ 
and discriminant $D_L$. We have $\xi^2=\alpha$, therefore $L$ is a subfield of $K$.

\vspace{1cm}

\begin{lemma}
$K$ is a composite of the totally real cubic number field $L$
and the imaginary quadratic field $M=\Q(i)$, where the discriminants of $L$ and $M$
are coprime.
\end{lemma}

\noindent
{\bf Proof.} 
Substituting $x=y-(n^2+5)/3$ into $g(x)$ and writing the resulting polynomial
in the form $h(y)=y^3+py+q$ we obtain
\[
D(h)=\left(\frac{q}{2}\right)^2+\left(\frac{p}{3}\right)^3=
-\frac{1}{108}(n^2+n+7)^2(n^2+n-1)^2<0,
\]  
whence the roots $\alpha_1,\alpha_2,\alpha_3$ of $g(x)$ are real numbers:
$L$ is a totally real cubic number field.

Set
\[
\zeta=(-n^2-2n^2-3n+3)\xi+(-n^3+n^2-5n+4)\xi^3+(-n+1)\xi^5.
\]
A Maple calculation implies
\[
\zeta^2=(-1)(n^2+n-1)^2,
\]
whence $\zeta=i\cdot (n^2+n-1)$. Hence $M=\Q(i)$ is also a subfield of $K$. 

An integral basis of $M$ is $(1,i)$, the discriminant of $M$ is $D_M=-4$.
The discriminant of $g(x)$ is
\[
D(g)=(n^2+n+7)^2(n^2+n-1)^2,
\]
which is always an odd number. $D_L$ is a divisor of $D(g)$, hence $D_L$ is also odd.
These imply, that $(D_M,D_L)=1$.
\hfill $\Box$

\vspace{1cm}

As a consequence of $(D_M,D_L)=1$ we obtain:

\begin{lemma}
If 
\[
(1,\beta_2,\beta_3)
\]
 is an integral basis of $L$ then
\[
(1,\beta_2,\beta_3,i,i\beta_2,i\beta_3)
\]
 is an integral basis of $K$.
\end{lemma}

\vspace{1cm}

The properties of $K$ allow to apply the method described in \cite{totrelext}.

Denote by $\alpha^{(j)}\; (j=1,2,3)$ the conjugates of $\alpha$, and let
$\beta_2^{(j)}, \beta_3^{(j)}$ be the conjugates of $\beta_2,\beta_3$, respectively,
corresponding to $\alpha^{(j)}$.

Let 
\begin{equation}
\gamma=x_1+x_2\beta_2+x_3\beta_3+iy_1+iy_2\beta_2+iy_3\beta_3\in\Z_K
\label{gamma}
\end{equation}
be arbitrary with $x_1,x_2,x_3,y_1,y_2,y_3\in\Z$.
Then
\[
\gamma^{(1,j)}=x_1+x_2\beta_2^{(j)}+x_3\beta_3^{(j)}
+iy_1+iy_2\beta_2^{(j)}+iy_3\beta_3^{(j)}
\]
and
\[
\gamma^{(2,j)}=x_1+x_2\beta_2^{(j)}+x_3\beta_3^{(j)}
-iy_1-iy_2\beta_2^{(j)}-iy_3\beta_3^{(j)}
\]
($j=1,2,3$) are the conjugates of $\gamma$.

\section{Auxiliary results}

We formulate the following general results of \cite{grs}
for our special case of a sextic field
$K$ being a composite of a totally real cubic field $L$
and an imaginary quadratic subfield $M=\Q(i)$.

\begin{lemma} For  $\gamma\in\Z_K$ generating $K$ over $\Q$ we have
\[
I(\gamma)=I_{K/M}(\gamma)\cdot J(\gamma)
\]
where
\[
I_{K/M}(\gamma)=(\Z_K:\Z_M[\gamma])=
\frac{1}{\sqrt{|N_{M/\Q}(D_{K/M})|}}
\prod_{i=1}^2\prod_{1\le j_1< j_2\le 3}|\gamma^{(i,j_1)}-\gamma^{(i,j_2)}|
\]
is the relative index of $\gamma$ and
\[
J(\gamma)=\frac{1}{|D_M|^{3/2}}
\prod_{j_1=1}^3\prod_{j_2=1}^3 |\gamma^{(1,j_1)}-\gamma^{(2,j_2)}|.
\]
\label{relindex}
\end{lemma}

If $\gamma$ generates a power integral basis in $K$, 
that is $I(\gamma)=1$, then $I_{K/M}(\gamma)=1$ and $J(\gamma)=1$.
Let $I_L(x_2,x_3)\in\Z[x,y]$ be the index form corresponding to the 
integral basis $(1,\beta_2,\beta_3)$ of $L$.

In \cite{totrelext} we showed that in our spacial case: 

\begin{lemma}
If $I_{K/M}(\gamma)=1$ then
\begin{equation}
|I_L(x_2,x_2)|\le 1,\;\; |I_L(y_2,y_3)|\le 1.
\label{thue}
\end{equation}
\end{lemma}

This is the main power of the method described in \cite{totrelext}.
The relative index form equation $I_{K/M}(\gamma)=1$ 
(in our case a cubic relative Thue equation over the quadratic subfield $M$,
cf. \cite{book}) implies absolute index form equations in $L$,
which are cubic Thue equations. Moreover, in the special case of our 
number field $K$, $J(\gamma)$ also factorizes:

\begin{lemma}
If $J(\gamma)=1$ then
\begin{equation}
N_{L/Q}(y_1+\beta_2 y_2+\beta_3 y_3)=\pm 1
\label{norm}
\end{equation}
and
\begin{equation}
P(\gamma)=\prod_{
\scriptsize{
\begin{array}{c}1\le j_1,j_2\le 3\\j_1\ne  j_2\end{array}
}
}|\gamma^{(1,j_1)}-\gamma^{(2,j_2)}|
=\pm 1.
\label{jeq}
\end{equation}
\end{lemma}

\noindent
{\bf Proof.}
In our case $|D_M|=4$ and
\[
\frac{1}{2^3}\prod_{j=1}^3 (\gamma^{(1,j)}-\gamma^{(2,j)})
\]
is a symmetric polynomial with rational integer coefficients,
equal to $N_{L/Q}(y_1+\beta_2 y_2+\beta_3 y_3)$. The corresponding 
factor of $J(\gamma)$ is $P(\gamma)$, also having rational integer coefficients. 
\hfill $\Box$

\section{The algorithm}

In view of the above statements, in order to determine all 
non-equivalent generators of power integral bases of $K$, we
perform the following steps for each parameter value $n$:

\begin{enumerate}
\item \label{a1}
Calculate an integer basis $(1,\beta_2,\beta_3)$ of $L$.
\item \label{a2}
Solve $I_L(x_2,x_3)=\pm 1$. Let $H$ be the set of solutions $(x_2,x_3)$.
\item \label{a3}
Let $H_0=H\cup \{(0,0)\}$.
\item \label{a4}
For all $(y_1,y_2)\in H_0$ calculate the corresponding $y_1$.
Let $H_1$ be the set of possible triples $(y_1,y_2,y_3)$.
\item \label{a5}
For all $(x_1,x_2)\in H_0$ and for all $(y_1,y_2,y_3)\in H_1$
construct $\gamma$ (cf. (\ref{gamma})) 
and test if $I_{K/M}(\gamma)=1$ and  $P(\gamma)=1$ hold.
\end{enumerate}

We performed all calculations in Maple, also the integral basis in Step \ref{a1}
was calculated by Maple.

In Step \ref{a2} we calculated the solutions $(x_2,x_3)$ of 
\begin{equation}
I_L(x_2,x_3)=\pm 1,\;\; x_2,x_3\in \Z,\; {\rm with} \;|x_2|,|x_3|\le 10^{100}.
\label{i23}
\end{equation} 
For larger parameters the coefficients of $I_L$ become extremely large,
therefore the complete resolution of (\ref{i23}) would have been very time consuming.

In Step \ref{a4}, for given $(y_2,y_3)$ we calculated $y_1$ using
equation (\ref{norm}) which is then a polynomial equation in $y_1$
with integer coefficients.

Step \ref{a5} is necessary to select the solutions from the set of 
possible solutions.

\vspace{1cm}

As an example for equation (\ref{i23}), for $n=140$ we provide here 
the integral basis \\ $(1,x,(x^2+58941x+118925)/138173)$ of $L$
and the corresponding index form equation:
\[
I(x_2,x_3)= -138173 x_2^3  - 137613 x_2^2  x3 - 44758 x_2 x_3^2  - 4777 x_3^3.
\]

\vspace{1cm}

\noindent
We performed two series of explicit numerical calculations.

\vspace{1cm}

\noindent
{\bf A.} $-100\le n\le 100$.\\
Calculating the solutions of (\ref{i23}) with $-100\le n\le 100$ 
took about 30 minutes,
out of which the calculation for the interval $-50\le n\le 50$
took only 1.5 minutes. 
This shows how the large coefficients slow down the calculations.

\noindent
{\bf B.} $n\in S$, where
\[
S=\{n\; : \; n\in [-1000,-100)\cup (100,1000],
 \;\; n^2+n+7 \; {\rm square}\; {\rm free}    \}.
\]
The set $S$ contains 1110 parameters $n$.
The reason to consider this set is that for all $n\in S$
we have the same type of integer basis.  Hence we can 
write equation (\ref{i23}) in a parametric form 
and we can perform also Step \ref{a4} and Step \ref{a5} in a parametric form.
It took 39 minutes to find the solution of $I(x_2,x_3)=\pm 1$
with $|x_2|,|x_3|\le 10^{100}$ for all the 1110 parameters $n\in S$.

\section{Results}

For the set {\bf A} of parameters $n$ our explicit calculations imply:
\begin{theorem}
\label{th1}
$n=-56,-14,-7,-5,-2,-1,0,1,4,6,13$ are the only values of $n$ with
$-100\le n\le 100$, such that $K$ admits generators of power integral bases
with coefficients $\le 10^{100}$ in the integral basis.
\end{theorem}

This statement was proved by explicit numerical calculations, following 
the above Algorithm. All data of generators of power integral bases
in the monogenic fields are listed in Section \ref{table}.

For the set {\bf B} of parameters we only calculated the solutions of 
equation (\ref{i23}) using Maple, 
and the remaining calculations were made in a parametric
form.

\begin{theorem}
For $n\in S$ there are no generators $\gamma$ (cf. (\ref{gamma}))  of power 
integral bases of $K$ 
with coefficients $x_2,x_3,y_1,y_2,y_3\in \Z$
having absolute values $<10^{100}$. 
\end{theorem}

\noindent
{\bf Remark.} We can not exclude the existence of monogenic fields $K$
for $|n|>100$, but they are certainly not of type $S$.

\noindent
{\bf Proof}. 
Our explicit calculations showed that for all $n\in S$ an integral basis 
of $L$ is given by
\[
\left(1,\alpha,\beta\right),
\]
where
\[
\beta=\frac{\alpha^2+(n^2-n+3)\alpha+n^2}{n^2+n-1}.
\]
The index form of $L$, corresponding to this integral basis is
\[
I_L(x_2,x_3)=(n^2+n-1)x_2^2+(n^2-3n-1)x_2^2x_3+(2-2n)x_2x_3^2+x_3^3.
\]
Let
\begin{equation}
T:=\{(0,0),(\pm 1,\pm n),(0,\pm 1),(\pm 1,\pm (n-1))\}.
\label{tmo}
\end{equation}
It is easy to check that $(x_2,x_3)\in T$ are solutions of $|I(x_2,x_3)|\le 1$.
For $n\in S$ we did not find any further solutions of (\ref{i23}).

For $n\in S$,  $(0,0)\ne (y_2,y_3)\in T$ there exist
no corresponding $y_1$. This can be shown by explicit calculations
using symmetric polynomials. For example, for $y_2=1,y_3=n$ the left hand side
of (\ref{norm}) is 
\[
N=y_1^3 + (2n - 5)y_1^2 + (n^2 - 7n + 6)y_1 - 2n^2 + 4n - 1.
\]
We have 
\[
N+1=(ny_1 + y_1^2 - 2n - 3y_1)(y_1 + n - 2)
\]
\[
N-1=(y_1 + n - 1)(ny_1 + y_1^2 - 2n - 4y_1 + 2)
\]
The above second degree factors are non-zero for $n\in S$. There remains
$y_1=2-n,1-n$ to test.  For these triplets $(y_1,y_2,y_3)$ and for all
$(x_2,x_3)\in T$ we calculated $P(\gamma)$ using again symmetric polynomials.
These $P(\gamma)$ are polynomials of $n$ with integer coefficients, for which
neither $P(\gamma)+1$, nor $P(\gamma)-1$ has integer roots in $n$. The other 
possible non-zero pairs $(y_2,y_3)$ were considered similarly.

For $(y_2,y_2)=(0,0)$ we obviously have $y_1=\pm 1$. For
$(y_1,y_2,y_3)=(\pm 1,0,0)$ we again tested all $(x_2,x_3)\in T$
and found that neither $P(\gamma)+1$, nor $P(\gamma)-1$ has integer roots
in $n$.
\hfill $\Box$

\section{Table}
\label{table}

Here we list the values of $n$ with $-100\le n\le 100$ for which 
we found generators $\gamma$ of power integral bases of $K$
with coefficients $<10^{100}$ in absolute value in the integral bases.
We display $n$ and the integral basis $(1,\beta_1,\beta_2)$ 
of $L$. We display the coefficients $(x_2,x_3,y_1,y_2,y_3)$
of non-equivalent generators $\gamma$ (cf. (\ref{gamma})) of power integral bases 
with respect to the integral basis $(1,\beta_1,\beta_2,i,i\beta_1,i\beta_2)$.

\noindent
$n=-56$, integral basis of $L$: 
$\displaystyle{\left(1,\frac{x+4}{7},\frac{x^2 + 9353x + 18531}{150871}\right)}$\\
$(x_2,x_3,y_1,y_2,y_3)=
                    ( 0, 0, 9, -49, 170),
                     ( 0, 0, 10, -49, 170),
                    ( 0, 0, -4, 17, -59),$ \\ $
                    ( 0, 0, -3, 17, -59),
                    ( 0, 0, -5, 32, -111),
                    ( 0, 0, -4, 32, -111)$

\vspace{0.5cm}

\noindent
$n=-14$, integral basis of $L$: 
$\displaystyle{\left(1,\frac{x+1}{3},\frac{x^2 + 32x + 1282}{1629}\right)}$\\
$(x_2,x_3,y_1,y_2,y_3)=
                    ( 0, 0, 12, -5, -16),
                    ( 0, 0, 13, -5, -16),
                    ( 0, 0, 11, -4, -13),$ \\ $
                    ( 0, 0, 12, -4, -13),
                     ( 0, 0, -25, 9, 29),
                     ( 0, 0, -24, 9, 29)$

\vspace{0.5cm}

\noindent
$n=-7$, integral basis of $L$: 
$\displaystyle{\left(1,x,\frac{x^2 + 141x + 254}{287}\right)}$\\
$(x_2,x_3,y_1,y_2,y_3)=
                    ( 0, 0, -19, -7, 23),
                    ( 0, 0, -18, -7, 23),
                     ( 0, 0, 7, 3, -10),$ \\ $
                     ( 0, 0, 8, 3, -10),
                     ( 0, 0, 11, 4, -13),
                     ( 0, 0, 12, 4, -13)$

\vspace{0.5cm}

\noindent
$n=-5$, integral basis of $L$: 
$\displaystyle{\left(1,x,\frac{x^2 + 14x + 25}{57}\right)}$\\
$(x_2,x_3,y_1,y_2,y_3)=
                     ( 0, 0, 0, -2, -7),
                     ( 0, 0, 1, -2, -7),
                     ( 0, 0, 1, -1, -4),$ \\ $
                     ( 0, 0, 2, -1, -4),
                     ( 0, 0, -4, 3, 11),
                     ( 0, 0, -3, 3, 11)$

\vspace{0.5cm}

\noindent
$n=-2$, integral basis of $L$: $(1,x,x^2)$\\
$(x_2,x_3,y_1,y_2,y_3)=
                    ( 0, 0, -5, -17, -2),
                    ( 0, 0, -4, -17, -2),
                     ( 0, 0, -4, -9, -1),$ \\ $
                     ( 0, 0, -3, -9, -1),
                      ( 0, 0, 9, 26, 3),
                      ( 0, 0, 10, 26, 3)$

\vspace{0.5cm}

\noindent
$n=-1$, integral basis of $L$: $(1,x,x^2)$\\
$(x_2,x_3,y_1,y_2,y_3)=
                     ( 0, 0, -4, -6, -1),
                     ( 0, 0, -3, -6, -1),
                     ( 0, 0, -2, -5, -1),$ \\ $
                     ( 0, 0, -1, -5, -1),
                      ( 0, 0, 4, 11, 2),
                      ( 0, 0, 5, 11, 2)$

\vspace{0.5cm}

\noindent
$n=0$, integral basis of $L$: $(1,x,x^2)$\\
$(x_2,x_3,y_1,y_2,y_3)=
                     ( 0, 0, -1, -3, -1),
                     ( 0, 0, 0, -3, -1),
                   ( 0, 0, -2, -1, 0),$ \\ $
                     ( 0, 0, -1, -1, 0),
                      ( 0, 0, 2, 4, 1),
                      ( 0, 0, 3, 4, 1)$

\vspace{0.5cm}

\noindent
$n=1$, integral basis of $L$: $(1,x,x^2)$\\
$(x_2,x_3,y_1,y_2,y_3)=
                     ( 0, 0, -1, -3, -1),
                     ( 0, 0, 0, -3, -1),
                     ( 0, 0, -3, -1, 0),$ \\ $
                     ( 0, 0, -2, -1, 0),
                      ( 0, 0, 1, 4, 1),
                      ( 0, 0, 2, 4, 1)$

\vspace{0.5cm}

\noindent
$n=4$, integral basis of $L$: 
$\displaystyle{\left(1,x,\frac{x^2 + 53x + 16}{57}\right)}$\\
$(x_2,x_3,y_1,y_2,y_3)=
                      ( 0, 0, 0, -4, 7),
                      ( 0, 0, 1, -4, 7),
                     ( 0, 0, -1, 1, -2),$ \\ $
                      ( 0, 0, 0, 1, -2),
                      ( 0, 0, 1, 3, -5),
                      ( 0, 0, 2, 3, -5)$

\vspace{0.5cm}

\noindent
$n=6$, integral basis of $L$: 
$\displaystyle{\left(1,x,\frac{x^2 + 74x + 200}{287}\right)}$\\
$(x_2,x_3,y_1,y_2,y_3)=
                    ( 0, 0, -11, -2, 17),
                    ( 0, 0, -10, -2, 17),
                      ( 0, 0, 4, 1, -9),$ \\ $
                      ( 0, 0, 5, 1, -9),
                      ( 0, 0, 5, 1, -8),
                      ( 0, 0, 6, 1, -8)$

\vspace{0.5cm}

\noindent
$n=13$, integral basis of $L$: 
$\displaystyle{\left(1,x,\frac{x^2 + 521x + 169}{1629}\right)}$\\
$(x_2,x_3,y_1,y_2,y_3)=
                     ( 0, 0, 4, -16, 25),
                     ( 0, 0, 5, -16, 25),
                     ( 0, 0, -3, 7, -11),$ \\ $
                     ( 0, 0, -2, 7, -11),
                     ( 0, 0, -2, 9, -14),
                     ( 0, 0, -1, 9, -14)$

\end{document}